\documentclass[11pt, a4paper]{article}

\usepackage[a4paper, left=3cm, right=3cm, top=2.5cm]{geometry}
\usepackage[utf8]{inputenc}
\usepackage[T1]{fontenc}
\usepackage[english,francais]{babel} 
\usepackage{xcolor} 
\usepackage[colorlinks, linkcolor=black, pagebackref, urlcolor=gray, citecolor=blue]{hyperref} 
\usepackage{amsthm}
\usepackage{amssymb} 
\usepackage{amsmath} 
\usepackage{dsfont} 
\usepackage{enumitem} 

\newcounter{compt}
\newtheorem{theoreme}[compt]{Théorème}
\newtheorem*{theoremenn}{Théorème} 
\newtheorem{lemme}[compt]{Lemme}

\begin{document}
\title{Théorème d'Erd\H{o}s-Kac pour les translatés\\
d'entiers ayant $k$ facteurs premiers}
\author{Élie \bsc{Goudout}}
\date{}
\maketitle
\vspace{-.5cm}

\selectlanguage{english}

\begin{abstract}
Let $x\geqslant 3$. For $1\leqslant n\leqslant x$ an integer, let $\omega(n)$ be its number of distinct prime factors. We show that $\omega(n-1)$ satisfies an Erd\H{o}s-Kac type theorem whenever $\omega(n)=k$ where $1\leqslant k\ll\log\log x$, thus extending a result of Halberstam.
\end{abstract}

\selectlanguage{francais}

\section{Présentation des résultats}

Pour $k, n\geqslant 1$ deux entiers et $x\geqslant 1$ un réel, on note $\omega(n)$ le nombre de facteurs premiers distincts de $n$,
\[\mathcal{E}_k(x):=\{n\leqslant x\,:\quad\omega(n)=k\},\]
et $\pi_k(x):=\#\mathcal{E}_k(x)$. On s'intéresse à une variante du théorème suivant. On définit
\[\Phi(y):=\frac{1}{\sqrt{2\pi}}\int_{-\infty}^y\mathrm{e}^{-t^2/2}\mathrm{d}t.\hspace{1cm}(y\in\mathbb{R})\]
\begin{theoremenn}[Erd\H{o}s \& Kac~\cite{erdoskac} ; Rényi \& Tur\'an~\cite{renyituran}]
Uniformément pour $x\geqslant 3$ et $y\in\mathbb{R}$, on a
\[\#\left\{n\leqslant x\,:\quad\omega(n)\leqslant\log_2 x+y\sqrt{\log_2 x}\right\}=\Phi(y)+O\left(\frac{1}{\sqrt{\smash[b]{\log_2 x}}}\right).\]
\end{theoremenn}

Dans~\cite{repstatfouten}, Fouvry \& Tenenbaum étudient la répartition de $\omega(n-1)$ lorsque $n$ vérifie une condition de friabilité. On s'inspire de leur méthode pour étudier la répartition de $\omega(n-1)$ lorsque $n$ a un nombre fixé de facteurs premiers. Le cas particulier où $n=p$ est premier a été traité par Halberstam~\cite{halberstamerdoskac}. Pour $x\geqslant 3$, $y\in\mathbb{R}$ et $k\geqslant 1$ un entier, on définit
\[\pi_k(x, y):=\#\left\{n\in\mathcal{E}_k(x)\,:\quad\omega(n-1)\leqslant\log_2 x+y\sqrt{\log_2 x}\right\}.\]

\begin{theoreme}\label{tqhgvsjbk}
Soit $R>0$ fixé. Uniformément pour $x\geqslant 3$, $1\leqslant k\leqslant R\log_2 x$ et $y\in\mathbb{R}$, on a
\[\pi_k(x, y)=\pi_k(x)\left\{\Phi(y)+O\left(\frac{\log_3 x}{\sqrt{\smash[b]{\log_2 x}}}\right)\right\}.\]
\end{theoreme}

On démontre ce théorème par la méthode des fonctions caractéristiques et l'inégalité de Berry-Esseen. On obtient le même terme d'erreur que Fouvry \& Tenenbaum. Dans les deux cas, on néglige la contribution des grands facteurs premiers de $n-1$ pour des raisons techniques, liées à la répartition de l'ensemble étudié dans les grandes progressions arithmétiques. C'est pour cela qu'on n'obtient pas le terme d'erreur $O((\log_2 x)^{-1/2})$ espéré.

En utilisant l'idée de Selberg pour estimer $\pi_k(x)$, l'étude directe de la série génératrice adéquate nous permet aussi de démontrer le théorème suivant, qui localise le nombre de \og petits\fg\ facteurs premier de $n-1$, lorsque $n\in\mathcal{E}_k(x)$. Pour $k\geqslant 1$, $\ell\geqslant 0$ et $1\leqslant w\leqslant x$, on note
\begin{align*}
\omega(n, w)&:=\sum_{\substack{p\vert n \\ p\leqslant w}}1,\hspace{1cm}(n\geqslant 1)\\
\pi_{k, \ell}(x, w)&:=\#\left\{n\in\mathcal{E}_k(x)\,:\quad\omega(n-1, w)=\ell\right\}.
\end{align*}

\begin{theoreme}\label{tdshjfk}
Il existe une constante absolue $\eta>0$ vérifiant l'énoncé suivant. Soit $R\geqslant 1$ fixé quelconque. Uniformément pour $3\leqslant w\leqslant x^{\eta(\log_3 x)/(R^2\log_2 x)}$, $1\leqslant k\leqslant R\log_2 x$ et $0\leqslant\ell\leqslant R\log_2 w$, on a
\begin{equation}\label{rfdbazhjbhkcsdioukj}
\pi_{k, \ell}(x, w)=\pi_k(x)\frac{(\log_2 w)^{\ell}}{\ell !\log w}\left\{h_k\left(\frac{\ell}{\log_2 w}\right)+O\left(\frac{k}{(\log_2 x)^2}+\frac{\ell+1}{(\log_2 w)^2}\right)\right\},
\end{equation}
où $h_k$ est la fonction entière définie par
\begin{equation}\label{hefjzdsuh}
h_k(z):=\mathrm{e}^{\gamma(z-1)}\prod_{p\geqslant 2}\left(1+\frac{z-1}{p+(k-1)/\log_2 x-1}\right)\left(1-\frac{1}{p}\right)^{z-1}.\hspace{1cm}(z\in\mathbb{C})
\end{equation}
\end{theoreme}

Ce théorème est non trivial uniquement lorsque $w$ est suffisamment grand. Le cas des petites valeurs de $w$ peut être traité via un léger ajustement de la preuve, par l'utilisation de~(\ref{tghejfduhn}). On note que lorsque $k=1$ et $\ell=0$, le terme principal est nul, ce qui est normal puisque $2$ divise tous les $p-1$ pour $p\geqslant 3$. De manière général, lorsque $k/\log_2 x+\ell/\log_2 w$ est petit, le facteur local du système $(n, n-1)$ a un impact important.

\section{Étude de la série génératrice}

On étudie la série génératrice associée au membre de gauche de~(\ref{rfdbazhjbhkcsdioukj}), dans le but d'utiliser la méthode de Selberg. Cela nous permet, dans la dernière section, de démontrer les théorèmes~\ref{tqhgvsjbk} et~\ref{tdshjfk}.

\subsection{Moyenne d'une fonction multiplicative sur $\mathcal{E}_k(x)$}

Le premier lemme que l'on démontre se déduit simplement de la méthode de Selberg-Delange, telle qu'exposée dans~\cite{tenenbaum}. Étant donnés $x\geqslant 2$, $N\geqslant 0$ et $c_1, c_2>0$, on définit
\[R_N(x):=R_N(x ; c_1, c_2):=\mathrm{e}^{-c_1\sqrt{\log x}}+\left(\frac{c_2N+1}{\log x}\right)^{N+1}.\]

\begin{lemme}\label{ghjkfsdkhsdjlij}
Soit $0<\varepsilon<1$ et $A, R>0$ fixés. Il existe des constantes $c_1, c_2>0$, pouvant dépendre de $R$, vérifiant l'énoncé suivant. Uniformément pour $N\geqslant 0$, $\vert\kappa\vert\leqslant A$, $x\geqslant 3$, $1\leqslant k\leqslant \vert\kappa\vert R\log_2x$ et $f$ une fonction multiplicative vérifiant
\begin{enumerate}[label=(\roman*), leftmargin=.8cm, align=left, labelsep=*]
\item\label{aaa} $\displaystyle \sum_{p\geqslant 2}\frac{\vert f(p)-\kappa\vert}{p^{1-\varepsilon}}\leqslant A$,
\item\label{bbb} $\displaystyle \sum_{p\geqslant 2}\sum_{\nu\geqslant 2}\frac{\vert f(p^{\nu})\vert}{p^{\nu(1-\varepsilon)}}\leqslant A$,
\end{enumerate}
on a
\begin{equation}\label{tghejfduhn}
\sum_{n\in\mathcal{E}_k(x)}f(n)=\frac{x}{\log x}\Bigg\{\sum_{0\leqslant j\leqslant N}\frac{Q_{f, \kappa, j, k}(\kappa\log_2 x)}{(\log x)^j}+O\left(\frac{(\vert\kappa\vert\log_2 x)^k}{k!}R_N(x)\right)\Bigg\},
\end{equation}
où les polynômes $Q_{f, \kappa, j, k}$ sont explicites, de degré au plus $k-1$ et ne dépendent que de $\kappa, f, j$ et $k$. De plus, en notant $r:=(k-1)/(\kappa\log_2 x)$, on a
\begin{equation}\label{tghzhzejzeh}
\sum_{n\in\mathcal{E}_k(x)}f(n)=\frac{x}{\log x}\frac{(\kappa\log_2 x)^{k-1}}{(k-1)!}\Bigg\{\lambda_{f, \kappa}(r)-\frac{r\lambda_{f, \kappa}''(r)}{2\log_2 x}+O\left(\frac{k^2}{(\kappa\log_2 x)^4}\right)\Bigg\},
\end{equation}
où $\lambda_{f, \kappa}$ est la fonction entière définie par
\[\lambda_{f, \kappa}(z):=\frac{\kappa}{\Gamma(\kappa z+1)}\prod_{p\geqslant 2}\left(1+z\sum_{\nu\geqslant 1}\frac{f(p^{\nu})}{p^{\nu}}\right)\left(1-\frac{1}{p}\right)^{\kappa z}.\hspace{1cm}(z\in\mathbb{C})\]
\end{lemme}

\begin{proof}
On suppose donnés les paramètres de l'énoncé. Pour $z\in\mathbb{C}$, on introduit, lorsque cela a un sens, les séries de Dirichlet
\begin{align*}
F(s)&:=\sum_{n\geqslant 1}\frac{z^{\omega(n)}f(n)}{n^s},\\
G(s ; z)&:=F(s)\zeta(s)^{-\kappa z}=\prod_{p\geqslant 2}\left(1+z\sum_{\nu\geqslant 1}\frac{f(p^{\nu})}{p^{\nu s}}\right)\left(1-\frac{1}{p^s}\right)^{\kappa z}.
\end{align*}
En notant $b(n)$ les coefficients de la série $G(\cdot ; z)$, pour tous $p\geqslant 2$ et $\nu\geqslant 1$, on a
\[b(p^{\nu})={\nu-1-\kappa z \choose \nu}+z\sum_{0\leqslant\mu<\nu}f(p^{\nu-\mu}){\mu-1-\kappa z \choose \mu}.\]
En particulier, $b(p)=z(f(p)-\kappa)$. Lorsque $\vert z\vert\leqslant R$, en utilisant la majoration
\[{\mu-1-\kappa z \choose \mu}\ll(1+\mu)^{1+\vert\kappa z\vert},\]
avec $\varepsilon$ qui vérifie~\ref{aaa} et~\ref{bbb}, on obtient uniformément
\[\sum_{p\geqslant 2}\sum_{\nu\geqslant 1}\frac{\vert b(p^{\nu})\vert}{p^{\nu(1-\varepsilon/3)}}<+\infty.,\]

L'estimation~(\ref{tghejfduhn}) se déduit alors des théorèmes \MakeUppercase{\romannumeral 2}.1.3 et \MakeUppercase{\romannumeral 2}.5.2 de~\cite{tenenbaum}, puis d'une version très légèrement modifiée -- afin tenir compte du paramètre $\kappa$ -- de la preuve de~\cite[th.~\MakeUppercase{\romannumeral 2}.6.3]{tenenbaum}. On en déduit~(\ref{tghzhzejzeh}) avec $N=0$, en effectuant un développement asymptotique de
\[Q_{f, \kappa, 0, k}(X)=\sum_{m+\ell=k-1}\frac{1}{m!\ell!}\lambda_{f, \kappa}^{(m)}(0)X^{\ell}\]
à deux termes, conformément à la note de fin du chapitre~\MakeUppercase{\romannumeral 2}.6 de~\cite{tenenbaum}.
\end{proof}

\subsection{Calcul de la série génératrice}\label{tghjefuhjki}

Pour $1\leqslant w\leqslant x$ et $k\geqslant 1$, on définit
\begin{align}\label{hjfkd}
f_k(z)&:=\sum_{\substack{n\in\mathcal{E}_k(x)}}z^{\omega(n-1,w)}.\hspace{1cm}(z\in\mathbb{C})
\end{align}
On omet les variables $w$ et $x$ de la notation afin de l'alléger. On rappelle aussi la définition~(\ref{hefjzdsuh}) de $h_k$.

\begin{lemme}\label{tghjzefsdvh}
Il existe une constante $\eta>0$ absolue vérifiant l'énoncé suivant. Soit $R\geqslant1$ fixé. Uniformément pour $3\leqslant w\leqslant x$, $1\leqslant k\ll R\log_2 x$ et $\vert z\vert\leqslant R$, avec $r:=(k-1)/\log_2 x$ et $u:=(\log x)/\log w$, on a
\begin{equation}\label{tgjdsuhnj}
\begin{split}
f_k(z)=\pi_k(x)(\log w)^{z-1}\Bigg\{h_k(z)+\frac{r\xi(z)}{\log_2 x}+O\bigg(\frac{k^2}{(\log_2 x)^4}&+\frac{1}{\log w}\bigg)\Bigg\}\\
+O\big(&xu^{-\eta u}(\log x)^{2R}+x^{9/10}\big),
\end{split}
\end{equation}
pour une fonction entière $\xi$ explicite, pouvant dépendre des paramètres $w, x$ et $k$, admettant $1$ pour zéro et uniformément bornée pour $\vert z\vert\leqslant R$.
\end{lemme}

On note qu'il est possible d'obtenir un développement asymptotique beaucoup plus précis de $f_k(z)$, de manière analogue à~(\ref{tghejfduhn}). La présence du terme d'ordre $r/\log_2 x$ nous est utile dans la démonstration du Théorème~\ref{tqhgvsjbk}.

\begin{proof}
On suppose donnés les paramètres de l'énoncé et on introduit la fonction multiplicative $g_z$ définie par
\[g_z(n):=\mu(n)^2(z-1)^{\omega(n)}.\hspace{1cm}(n\geqslant 1)\]
Puisque l'on a $1\ast g_z=z^{\omega(\cdot)}$, on obtient
\begin{align*}
f_k(z)&=\sum_{n\in\mathcal{E}_k(x)}\sum_{\substack{q\vert n-1 \\ P^+(q)\leqslant w}}g_z(q)\\
&=\sum_{\substack{q\leqslant x^{1/3} \\ P^+(q)\leqslant w}}g_z(q)\sum_{\substack{n\in\mathcal{E}_k(x) \\ n\equiv 1[q]}}1+R_1,
\end{align*}
où l'on a posé 
\[R_1:=\sum_{n\in\mathcal{E}_k(x)}1\sum_{\substack{q>x^{1/3} \\ q\vert n-1 \\ P^+(q)\leqslant w}}g_z(q).\]
Pour $n\geqslant 1$, on note
\[n_{w}:=\prod_{\substack{p^{\nu}\Vert n \\ p\leqslant w}}p^{\nu}\]
la partie $w$-friable de $n$. Dans $R_1$, la deuxième somme est non vide seulement lorsque $(n-1)_w>x^{1/3}$, et dans ce cas, on a $(1\ast\vert g_z\vert)(n-1)\leqslant (R+2)^{\omega(n-1)}$. Ainsi,
\begin{align*}
R_1&\ll\sum_{\substack{n\leqslant x \\ n_{w}>x^{1/3}}}(R+2)^{\omega(n)}\ll\sum_{\substack{x^{1/3}<a\leqslant x \\ P^+(a)\leqslant w}}\sum_{\substack{b\leqslant x/a \\ P^-(b)>w}}(R+2)^{\omega(a)+\omega(b)}\\
&\ll \frac{x}{\log x}\left(\frac{\log x}{\log w}\right)^{R+2}\sum_{\substack{n\leqslant x \\ n_{w}>x^{1/3}}}\frac{(R+2)^{\omega(a)}}{a}.
\end{align*}
Or, avec l'astuce de Rankin, pour tout $\alpha\in[0, 1/2]$, on a
\[\sum_{\substack{n\leqslant x \\ n_{w}>x^{1/3}}}\frac{(R+2)^{\omega(a)}}{a}\ll x^{-\alpha/3}\exp\bigg((R+2)\sum_{p\leqslant w}\frac{1}{p^{1-\alpha}-1}\bigg).\]
Avec~\cite[lem.~2]{tenenbaumprobdivinter} et $\alpha=\min\big(\log(2u)/\log w, 1/2\big)$, il vient $R_1\ll xu^{-\eta u}(\log x)^{2R}+x^{9/10}$ pour une certaine constante absolue $\eta>0$. En posant
\begin{align*}
R_2&:=\sum_{\substack{q\leqslant x^{1/3} \\ P^+(q)\leqslant w}}g_z(q)\bigg(\sum_{\substack{n\in\mathcal{E}_k(x) \\ n\equiv 1[q]}}1-\frac{1}{\varphi(q)}\sum_{\substack{n\in\mathcal{E}_k(x) \\ (n,q)=1}}1\bigg),\\
R_3&:=\sum_{\substack{q> x^{1/3} \\ P^+(q)\leqslant w}}\frac{g_z(q)}{\varphi(q)}\sum_{\substack{n\in\mathcal{E}_k(x) \\ (n,q)=1}}1,
\end{align*}
on obtient
\begin{equation}\label{hejfkdsijk}
f_k(z)=\sum_{n\in\mathcal{E}_k(x)}\prod_{\substack{p\leqslant w \\ p\nmid n}}\left(1+\frac{z-1}{p-1}\right)+R_1+R_2-R_3.
\end{equation}
De manière analogue à $R_1$, avec l'astuce de Rankin on a $R_3\ll xu^{-\eta u}(\log x)^{2R}+x^{9/10}$. On majore $R_2$ avec~\cite[th.~1]{wolkezhan} (pour un résultat analogue sur $\Omega$, voir~\cite{timofeevkhripunova}). On obtient, pour tout $A>0$,
\begin{align*}
R_2&\ll_A\Bigg(\sum_{\substack{q\leqslant x^{1/3}}}\vert g_z(q)\vert^2\bigg\vert\sum_{\substack{n\in\mathcal{E}_k(x) \\ n\equiv 1[q]}}1-\frac{1}{\varphi(q)}\sum_{\substack{n\in\mathcal{E}_k(x) \\ (n,q)=1}}1\bigg\vert\Bigg)^{1/2}\frac{\pi_k(x)^{1/2}}{(\log x)^A}\\
&\ll_A\frac{x}{(\log x)^{A}}\bigg(\sum_{\substack{q\leqslant x^{1/3}}}\frac{(R+1)^{2\omega(q)}}{\varphi(q)}\bigg)^{1/2}\\
&\ll_A\frac{x}{(\log x)^{A-R-1}}\\
&\ll\pi_k(x)(\log w)^{\Re\mathrm{e}\,z-2}.
\end{align*}

On suppose dans un premier temps que $\vert p+z-2\vert>1/3$ pour tout $p\geqslant 2$. Afin d'estimer la somme de~(\ref{hejfkdsijk}), on considère la fonction multiplicative $f$ définie par
\begin{equation}\label{tghefjkzok}
f(n):=\prod_{\substack{p\leqslant w \\ p\vert n}}\frac{p-1}{p+z-2}.
\end{equation}
Puisqu'elle vérifie les hypothèses du Lemme~\ref{ghjkfsdkhsdjlij} pour $\kappa=1$, avec~(\ref{tghzhzejzeh}) on obtient
\begin{multline*}
\sum_{n\in\mathcal{E}_k(x)}\prod_{\substack{p\leqslant w \\ p\nmid n}}\left(1+\frac{z-1}{p-1}\right)=\frac{x}{\log x}\frac{(\log_2 x)^{k-1}}{(k-1)!}\prod_{p\leqslant w}\left(1+\frac{z-1}{p-1}\right)\\
\times\left\{\lambda_z(r)-\frac{r\lambda_z''(r)}{2\log_2 x}+O\left(\frac{k^2}{(\log_2 x)^4}\right)\right\},
\end{multline*}
où l'on a posé $r:=(k-1)/\log_2 x$ et
\[\lambda_z(r):=\prod_{p\leqslant w}\left(1+\frac{r}{p+z-2}\right)\left(1-\frac{1}{p}\right)^r\prod_{p>w}\left(1+\frac{r}{p-1}\right)\left(1-\frac{1}{p}\right)^r.\]
On pose $\emph{a priori}$
\[\xi(z):=\frac{h_k(z)\lambda_1''(r)}{2\lambda_1(r)}-\frac{\lambda_z''(r)}{2\lambda_1(r)}(\log w)^{1-z}\prod_{p\leqslant w}\left(1+\frac{z-1}{p-1}\right).\]
Cela définit une fonction entière, admettant $1$ pour zéro et bornée uniformément en $w, x$ et $k$ sous les conditions de l'énoncé. En effet, les pôles de $z\mapsto\lambda''_z(r)$ sont compensés par les zéros du produit eulérien et $\lambda_1(r)\asymp 1$ pour $0\leqslant r\leqslant R$. On en déduit l'estimation désirée pour $f_k(z)$ puisque l'on a d'une part,
\[\pi_k(x)=\frac{x}{\log x}\frac{(\log_2 x)^{k-1}}{(k-1)!}\left\{\lambda_1(r)-\frac{r\lambda''_1(r)}{2\log_2 x}+O\left(\frac{k^2}{(\log_2 x)^4}\right)\right\},\]
et d'autre part,
\[(\mathrm{e}^{\gamma}\log w)^{z-1}\prod_{p\leqslant w}\left(1-\frac{1}{p}\right)^{z-1}=1+O\left(\frac{1}{\log w}\right),\]
d'après la troisième formule de Mertens, uniformément pour $\vert z\vert\leqslant R$.

Pour compléter la preuve, on traite le cas où il existe un $p_0$ tel que $\vert p_0+z-2\vert\leqslant 1/3$. Puisque la définition~(\ref{tghefjkzok}) n'est plus nécessairement valide, on traite spécifiquement le facteur eulérien de $p_0$ dans~(\ref{hejfkdsijk}). Pour cela, on écrit
\[\prod_{\substack{p\leqslant w \\ p\nmid n}}\left(1+\frac{z-1}{p-1}\right)=\prod_{\substack{p\leqslant w \\ p\neq p_0}}\left(1+\frac{z-1}{p-1}\right)\Bigg\{\prod_{\substack{p\leqslant w \\ p\vert n \\ p\neq p_0}}\frac{p-1}{p+z-2}+\mathds{1}_{p_0\nmid n}\frac{z-1}{p_0-1}\prod_{\substack{p\leqslant w \\ p\vert n}}\frac{p-1}{p+z-2}\Bigg\},\]
qui est une combinaison linéaire de deux fonctions multiplicatives vérifiant les hypothèses du Lemme~\ref{ghjkfsdkhsdjlij}, uniformément en $p_0$ et $\vert z\vert\leqslant R$. La même utilisation du Lemme~\ref{ghjkfsdkhsdjlij} que précédemment permet alors de conclure par un calcul analogue.
\end{proof}

\section{Lois locales et répartition}\label{tshflsd}

On démontre le Théorème~\ref{tdshjfk}, puis le Théorème~\ref{tqhgvsjbk}.
\begin{proof}[Démonstration du Théorème~\ref{tdshjfk}]
On suppose donnés les paramètres de l'énoncé. La quantité $\pi_{k, \ell}(x, w)$ correspond au coefficient de $x^{\ell}$ dans $f_k(z)$, qui est
\[\pi_{k, \ell}(x, w)=\frac{1}{2i\pi}\oint_{\vert z\vert=\rho}\frac{f_k(z)}{z^{\ell+1}}\mathrm{d}z.\]
pour tout $\rho>0$. Pour la démonstration, on majore
\[\frac{r\xi(z)}{\log_2 x}\ll\frac{k}{(\log_2 x)^2}.\]
En posant
\[Q(X):=\sum_{a+b=\ell}\frac{h_k^{(a)}(0)X^b}{a!b!},\]
le terme principal vaut
\[\frac{\pi_k(x)}{\log w}Q(\log_2 w).\]
Comme il est détaillé dans la preuve de~\cite[th.~\MakeUppercase{\romannumeral 2}.6.3(6.13)]{tenenbaum}, puisque $\vert h_k''(z)\vert\ll 1$ pour $\vert z\vert\leqslant R$ uniformément en $k$, on a
\[Q(\log_2 x)=\frac{(\log_2 w)^{\ell}}{\ell !}\left\{h_k\left(\frac{\ell}{\log_2 w}\right)+O\left(\frac{\ell}{(\log_2 w)^2}\right)\right\}.\]
Il suffit donc de traiter le terme d'erreur de~(\ref{tgjdsuhnj}). Le cas $\ell=0$ est trivial. Lorsque $\ell\geqslant 1$, d'après la majoration de~\cite[\MakeUppercase{\romannumeral 2}.(6.14)]{tenenbaum}, avec $\rho:=\ell/\log_2 w$, on a
\[\oint_{\vert z\vert=\rho}(\log w)^{\Re\mathrm{e}\,z}\vert z\vert^{-\ell-1}\vert\mathrm{d}z\vert\ll\frac{(\log_2 w)^{\ell}}{\ell !}.\]
Pour conclure, puisque $R\geqslant 1$ et $1\leqslant k\leqslant R\log_2 x$, il suffit de remarquer que
\[\pi_k(x)\frac{(\log_2 w)^{\ell}}{\ell !\log w}\left(\frac{k}{(\log_2 x)^2}+\frac{\ell+1}{(\log_2w)^2}\right)\gg\frac{x}{(\log x)^{10R^2}}\]
dès que $\pi_k(x)>0$.
\end{proof}

\begin{proof}[Démonstration du théorème~\ref{tqhgvsjbk}]
On suppose donnés les paramètres de l'énoncé. On suit essentiellement la démonstration de~\cite[cor.~5]{repstatfouten}, en adaptant le raisonnement à notre problème. Pour toute la démonstration, on choisit
\begin{equation}\label{tgefzhdjsk}
w:=\exp\left(\frac{\log x}{(\log_2 x)^2}\right),
\end{equation}
et on se donne une constante $C>0$, à fixer plus tard. Pour $x\geqslant 2$ et $y\in\mathbb{R}$, on pose
\begin{align*}
\tilde{\pi}_k(x, y)&:=\#\left\{n\in\mathcal{E}_k(x)\,:\quad\omega(n-1, w)\leqslant \log_2 x+y\sqrt{\log_2 x}\right\},\\
D_k(x)&:=\#\left\{n\in\mathcal{E}_k(x)\,:\quad\omega(n-1)-\omega(n-1, w)>C\log_3 x\right\}.
\end{align*}
On montre dans un premier temps que $\tilde{\pi}_k(x, y)$ est une bonne approximation de $\pi_k(x, y)$ lorsque $C$ est suffisamment grand. Pour cela, on observe que l'on a,
\begin{equation}\label{tgzhdjiu}
\tilde{\pi}_k\Big(x, y-\frac{C\log_3 x}{\sqrt{\smash[b]{\log_2 x}}}\Big)-D_k(x)\leqslant \pi_k(x, y)\leqslant\tilde{\pi}_k(x, y),
\end{equation}
et on majore $D_k$. On utilise l'astuce de Rankin, afin de majorer, pour tout $n\leqslant x$,
\begin{align*}
\mathds{1}_{\omega(n-1)-\omega(n-1, w)>C\log_3 x}&\leqslant 2^{\omega(n-1)-\omega(n-1, w)-C\log_3 x}\\
&\leqslant(\log_2 x)^{-C\log 2}\tau_2\big(\prod_{\substack{p^{\nu}\Vert n \\ p>w}}p^{\nu}\big).
\end{align*}
D'après~\cite{Landreau}, pour tous $\delta>0$ et $n\geqslant 1$, on a
\[\tau_2(n)\ll_{\delta}\sum_{\substack{q\vert n \\ q\leqslant n^{\delta}}}\tau_2(q)^{1/\delta}.\]
Avec $\delta=1/3$, en posant
\[R_4:=(\log_2 x)^{-C\log 2}\sum_{\substack{q\leqslant x^{1/3} \\ P^-(q)>w}}\tau_2(q)^{3}\bigg(\sum_{\substack{n\in\mathcal{E}_k(x) \\ n\equiv 1[q]}}1-\frac{1}{\varphi(q)}\sum_{\substack{n\in\mathcal{E}_k(x) \\ (n, q)=1}}1\bigg),\]
on obtient alors
\begin{align}
D_k(x)&\ll(\log_2 x)^{-C\log 2}\sum_{\substack{q\leqslant x^{1/3} \\ P^-(q)>w}}\tau_2(q)^{3}\sum_{\substack{n\in\mathcal{E}_k(x) \\ n\equiv 1[q]}}1\nonumber\\
&\ll(\log_2 x)^{-C\log 2}\sum_{n\in\mathcal{E}_k(x)}\sum_{\substack{q\leqslant x^{1/3} \\ (q, n)=1 \\ P^-(q)>w}}\frac{\tau_2(q)^{3}}{\varphi(q)}+R_4\nonumber\\
&\ll\frac{\pi_k(x)}{(\log_2 x)^{C\log 2}}\left(\frac{\log x}{\log w}\right)^{8}+R_4\nonumber\\
&\ll\frac{\pi_k(x)}{\log_2 x},\label{tghfjdshn}
\end{align}
lorsque $C$ est suffisamment grand, où l'on a majoré $R_4$ de manière analogue à $R_2$ dans la section~\ref{tghjefuhjki}.

Il suffit donc d'estimer $\tilde{\pi}_k(x, y)$. On utilise pour cela la méthode des fonctions caractéristiques. Pour des raisons techniques, en notant
\[\tilde{\pi}_k^{\star}(x, y):=\#\left\{n\in\mathcal{E}_k(x)\,:\quad\omega(n-1, w)\leqslant \log_2 w+y\sqrt{\log_2 w}\right\},\]
on montre en fait
\[\sup_{y\in\mathbb{R}}\vert\tilde{\pi}_k^{\star}(x, y)-\pi_k(x)\Phi(y)\vert\ll\frac{1}{\sqrt{\smash[b]{\log_2 x}}}.\]
Cela est suffisant puisque $\log_2 w=\log_2 x+O(\log_3 x)$ avec~(\ref{tgefzhdjsk}). D'après l'inégalité de Berry-Esseen telle qu'énoncée dans~\cite{tenenbaum}, en posant $T:=\log_2 w$, avec la notation~(\ref{hjfkd}) on a
\begin{equation}\label{tghfsjdk}
\sup_{y\in\mathbb{R}}\vert\tilde{\pi}_k^{\star}(x, y)-\pi_k(x)\Phi(y)\vert\ll\frac{\pi_k(x)}{\sqrt{T}}+\int_{-\sqrt{T}}^{\sqrt{T}}\left\vert \mathrm{e}^{-it\sqrt{T}}f_k\big(\mathrm{e}^{it/\sqrt{T}}\big)-\pi_k(x)\mathrm{e}^{-t^2/2}\right\vert\frac{\mathrm{d}t}{\vert t\vert}.
\end{equation}
Pour $\vartheta\in\mathbb{R}$, on a $\mathrm{e}^{i\vartheta}=1+O(\vartheta)$. Ainsi, lorsque $\vert t\vert\leqslant1/\log x$,
\begin{align*}
\left\vert \mathrm{e}^{-it\sqrt{T}}f_k\big(\mathrm{e}^{it/\sqrt{T}}\big)-\pi_k(x)\mathrm{e}^{-t^2/2}\right\vert\ll&\frac{\vert t\vert}{\sqrt{T}}\sum_{n\in\mathcal{E}_k(x)}\vert\omega(n-1, w)-T\vert+\vert t\vert^2\pi_k(x)\\
\ll&\vert t\vert \pi_k(x)\left(\frac{\log w}{\sqrt{T}}+\sqrt{T}+\vert t\vert\right),
\end{align*}
où l'on a majoré trivialement $\omega(n-1, w)\ll\log w$ . On a donc
\begin{equation}\label{tgehzj}
\int_{-1/\log x}^{1/\log x}\left\vert \mathrm{e}^{-it\sqrt{T}}f_k\big(\mathrm{e}^{it/\sqrt{T}}\big)-\pi_k(x)\mathrm{e}^{-t^2/2}\right\vert\frac{\mathrm{d}t}{\vert t\vert}\ll\frac{\pi_k(x)}{\sqrt{T}}.
\end{equation}
Pour $1/\log x\leqslant\vert t\vert\leqslant T^{1/3}$, avec le Lemme~\ref{tghjzefsdvh}, puisque $h_k(1)=1$ et $\xi(1)=0$, on a
\[f_k\big(\mathrm{e}^{it/\sqrt{T}}\big)=\pi_k(x)\mathrm{e}^{it\sqrt{T}}\mathrm{e}^{-t^2/2}\left\{1+O\left(\frac{\vert t\vert +\vert t\vert^3}{\sqrt{T}}+\frac{1}{(\log_2 x)^2}\right)\right\}.\]
Par ailleurs, en vertu de l'inégalité $\cos(\vartheta)-1\leqslant-2\vartheta^2/\pi^2$, valable pour $-1\leqslant\vartheta\leqslant 1$, lorsque $T^{1/3}<\vert t\vert\leqslant\sqrt{T}$, le Lemme~\ref{tghjzefsdvh} fournit la majoration
\[f_k\big(\mathrm{e}^{it/\sqrt{T}}\big)\ll\pi_k(x)\mathrm{e}^{-2t^2/\pi^2}.\]
On obtient donc
\begin{equation}\label{tgefhzjdsuyh}
\int_{1/\log x<\vert t\vert\leqslant\sqrt{T}}\left\vert \mathrm{e}^{-it\sqrt{T}}f_k\big(\mathrm{e}^{it/\sqrt{T}}\big)-\pi_k(x)\mathrm{e}^{-t^2/2}\right\vert\frac{\mathrm{d}t}{\vert t\vert}\ll\frac{\pi_k(x)}{\sqrt{T}}.
\end{equation}
Le résultat découle finalement des inégalités~(\ref{tgzhdjiu}) à~(\ref{tgefhzjdsuyh}). Il est possible d'obtenir un développement asymptotique plus précis pour $\pi_k^{\ast}(x, y)$. Malheureusement, cela ne permet pas d'améliorer le Théorème~\ref{tqhgvsjbk}.
\end{proof}

\noindent\textbf{Remerciements.} Je remercie Sary Drappeau et Berke Topacogullari pour m'avoir présenté ce sujet et pour nos discussions. Je remercie également Régis de la Bretèche pour ses précieuses suggestions.

\bibliographystyle{alpha-fr}
\hypersetup{linkcolor=gray}
\newcommand{\mapolicebackref}[1]{%
         \hspace*{-5pt}{\textcolor{gray}{\small$\uparrow$~#1}}
}
\renewcommand*{\backref}[1]{
\mapolicebackref{#1}
}

\bibliography{Bibliographie}

\noindent\bsc{Institut de Math\'ematiques de Jussieu-PRG, Universit\'e Paris Diderot,
Sorbonne Paris Cit\'e, 75013 Paris, France}\\

\noindent\textit{E-mail :} \url{eliegoudout@hotmail.com}

\end{document}